# ON SURROGATE LOSS FUNCTIONS AND $F$-DIVERGENCES[1]

By XuanLong Nguyen, Martin J. Wainwright
and Michael I. Jordan

*Duke University and Statistical and Applied Mathematical Sciences Institute (SAMSI), University of California, Berkeley and University of California, Berkeley*

The goal of binary classification is to estimate a discriminant function $\gamma$ from observations of covariate vectors and corresponding binary labels. We consider an elaboration of this problem in which the covariates are not available directly but are transformed by a dimensionality-reducing quantizer $Q$. We present conditions on loss functions such that empirical risk minimization yields Bayes consistency when both the discriminant function and the quantizer are estimated. These conditions are stated in terms of a general correspondence between loss functions and a class of functionals known as Ali-Silvey or $f$-divergence functionals. Whereas this correspondence was established by Blackwell [*Proc. 2nd Berkeley Symp. Probab. Statist.* **1** (1951) 93–102. Univ. California Press, Berkeley] for the 0–1 loss, we extend the correspondence to the broader class of surrogate loss functions that play a key role in the general theory of Bayes consistency for binary classification. Our result makes it possible to pick out the (strict) subset of surrogate loss functions that yield Bayes consistency for joint estimation of the discriminant function and the quantizer.

**1. Introduction.** Consider the classical problem of binary classification: given a pair of random variables $(X, Y) \in (\mathcal{X}, \mathcal{Y})$, where $\mathcal{X}$ is a Borel subset of $\mathbb{R}^d$ and $\mathcal{Y} = \{-1, +1\}$, and given of a set of samples $\{(X_1, Y_1), \ldots, (X_n, Y_n)\}$, the goal is to estimate a discriminant function that predicts the binary label $Y$ given the covariate vector $X$. The accuracy of any discriminant function

Received September 2007; revised September 2007.
[1]Supported in part by grants from Intel Corporation, Microsoft Research, Grant 0412995 from the National Science Foundation, and an Alfred P. Sloan Foundation Fellowship (MJW).
*AMS 2000 subject classifications.* 62G10, 68Q32, 62K05.
*Key words and phrases.* Binary classification, discriminant analysis, surrogate losses, $f$-divergences, Ali-Silvey divergences, quantizer design, nonparametric decentralized detection, statistical machine learning, Bayes consistency.

is generally assessed in terms of 0–1 loss as follows. Letting $\mathbb{P}$ denote the distribution of $(X,Y)$, and letting $\gamma : \mathcal{X} \to \mathbb{R}$ denote a given discriminant function, we seek to minimize the expectation of the 0–1 loss; that is, the error probability $\mathbb{P}(Y \neq \operatorname{sign}(\gamma(X)))$.[2] Unfortunately, the 0–1 loss is a nonconvex function, and practical classification algorithms, such as boosting and the support vector machine, are based on relaxing the 0–1 loss to a convex upper bound or approximation, yielding a *surrogate loss function* to which empirical risk minimization procedures can be applied. A significant achievement of the recent literature on binary classification has been the delineation of necessary and sufficient conditions under which such relaxations yield Bayes consistency [2, 9, 12, 13, 19, 22].

In many practical applications, this classical formulation of binary classification is elaborated to include an additional stage of "feature selection" or "dimension reduction," in which the covariate vector $X$ is transformed into a vector $Z$ according to a data-dependent mapping $Q$. An interesting example of this more elaborate formulation is a "distributed detection" problem, in which individual components of the $d$-dimensional covariate vector are measured at spatially separated locations, and there are communication constraints that limit the rate at which the measurements can be forwarded to a central location where the classification decision is made [21]. This communication-constrained setting imposes severe constraints on the choice of $Q$: any mapping $Q$ must be a separable function, specified by a collection of $d$ univariate, discrete-valued functions that are applied componentwise to $X$. The goal of decentralized detection is to specify and analyze data-dependent procedures for choosing such functions, which are typically referred to as "quantizers." More generally, we may abstract the essential ingredients of this problem and consider a problem of experimental design, in which $Q$ is taken to be a possibly stochastic mapping $\mathcal{X} \to \mathcal{Z}$, chosen from some constrained class $\mathcal{Q}$ of possible quantizers. In this setting, the discriminant function is a mapping $\gamma : \mathcal{Z} \to \mathbb{R}$, chosen from the class $\Gamma$ of all measurable functions on $\mathcal{Z}$. Overall, the problem is to simultaneously determine both the mapping $Q$ and the discriminant function $\gamma$, using the data $\{(X_1,Y_1),\ldots,(X_n,Y_n)\}$, so as to jointly minimize the Bayes error $R_{\text{Bayes}}(\gamma, Q) := \mathbb{P}(Y \neq \operatorname{sign}(\gamma(Z)))$.

As alluded to above, when $Q$ is fixed, it is possible to give general conditions under which relaxations of 0–1 loss yield Bayes consistency. As we will show in the current paper, however, these conditions *no longer suffice* to yield consistency in the more general setting, in which the choice of $Q$ is also optimized. Rather, in the setting of jointly estimating the discriminant function $\gamma$ and optimizing the quantizer $Q$, new conditions need to be imposed. It is the goal of the current paper to present such conditions and, moreover, to provide a general theoretical understanding of their origin. Such an understanding turns out to repose not only on analytic properties of surrogate loss functions (as in the $Q$-fixed case), but on a relationship between the family of surrogate loss functions and another class of functions known as $f$-*divergences* [1, 7]. In rough terms, an $f$-divergence between two distributions is defined by the expectation of a convex function of their likelihood ratio. Examples include the Hellinger distance, the total variational distance, Kullback–Leibler divergence and Chernoff distance, as well as various other divergences popular in the information theory literature [20]. In our setting, these $f$-divergences are applied to the class-conditional distributions induced by applying a fixed quantizer $Q$.

---

[2]We use the convention that $\operatorname{sign}(\alpha) = 1$ if $\alpha > 0$ and $-1$ otherwise.

An early hint of the relationship between surrogate losses and $f$-divergences can be found in a seminal paper of Blackwell [3]. In our language, Blackwell's result can be stated in the following way: if a quantizer $Q_A$ induces class-conditional distributions whose $f$-divergence is greater than the $f$-divergence induced by a quantizer $Q_B$, then there exists some set of prior probabilities for the class labels such that $Q_A$ results in a smaller probability of error than $Q_B$. This result suggests that any analysis of quantization procedures based on 0–1 and surrogate loss functions might usefully attempt to relate surrogate loss functions to $f$-divergences. Our analysis shows that this is indeed a fruitful suggestion, and that Blackwell's idea takes its most powerful form when we move beyond 0–1 loss to consider the full set of surrogate loss functions studied in the recent binary classification literature.

Blackwell's result [3] has had significant historical impact on the signal processing literature (and thence on the distributed detection literature). Consider, in a manner complementary to the standard binary classification setting in which the quantizer $Q$ is assumed known, the setting in which the discriminant function $\gamma$ is assumed known and only the quantizer $Q$ is to be estimated. This is a standard problem in the signal processing literature (see, e.g., [10, 11, 17]), and solution strategies typically involve the selection of a specific $f$-divergence to be optimized. Typically, the choice of an $f$-divergence is made somewhat heuristically, based on the grounds of analytic convenience, computational convenience or asymptotic arguments.

Our results in effect provide a broader and more rigorous framework for justifying the use of various $f$-divergences in solving quantizer design problems. We broaden the problem to consider the joint estimation of the discriminant function and the quantizer. We adopt a decision-theoretic perspective in which we aim to minimize the expectation of 0–1 loss, but we relax to surrogate loss functions that are convex approximations of 0–1 loss, with the goal of obtaining computationally tractable minimization procedures. By relating the family of surrogate loss functions to the family of $f$-divergences, we are able to specify equivalence classes of surrogate loss functions. The conditions that we present for Bayes consistency are expressed in terms of these equivalence classes.

1.1. *Our contributions.* In order to state our contributions more precisely, let us introduce some notation and definitions. Given the distribution $\mathbb{P}$ of the pair $(X, Y)$, consider a discrete space $\mathcal{Z}$, and let $Q(z|x)$ denote a *quantizer*—a conditional probability distribution on $\mathcal{Z}$ for almost all $x$. Let $\mu$ and $\pi$ denote measures over $\mathcal{Z}$ that are induced by $Q$ as follows:

$$\mu(z) := \mathbb{P}(Y = 1, Z = z) = p \int_x Q(z|x) \, d\mathbb{P}(x|Y = 1), \tag{1a}$$

$$\pi(z) := \mathbb{P}(Y = -1, Z = z) = q \int_x Q(z|x) \, d\mathbb{P}(x|Y = -1), \tag{1b}$$

where $p$ and $q$ denote the prior probabilities $p = \mathbb{P}(Y = 1)$ and $q = \mathbb{P}(Y = -1)$. We assume that $Q$ is restricted to some constrained class $\mathcal{Q}$, such that both $\mu$ and $\pi$ are strictly positive measures.

An $f$-divergence is defined as

$$I_f(\mu, \pi) := \sum_z \pi(z) f\left(\frac{\mu(z)}{\pi(z)}\right), \tag{2}$$

where $f : [0, +\infty) \to \mathbb{R} \cup \{+\infty\}$ is a continuous convex function. Different choices of convex $f$ lead to different divergence functionals [1, 7].

The loss functions that we consider are known as *margin-based* loss functions. Specifically, we study convex loss functions $\phi(y, \gamma(z))$ that are of the form $\phi(y\gamma(z))$, where the product $y\gamma(z)$ is known as the *margin*. Note in particular that 0–1 loss can be written in this form, since $\phi_{0-1}(y, \gamma(z)) = \mathbb{I}(y\gamma(z) \leq 0)$. Given such a margin-based loss function, we define the $\phi$-*risk* $R_\phi(\gamma, Q) = \mathbb{E}\phi(Y\gamma(Z))$. Statistical procedures will be defined in terms of minimizers of $R_\phi$ with respect to the arguments $\gamma$ and $Q$, with the expectation replaced by an empirical expectation defined by samples $\{(X_1, Y_1), \ldots, (X_n, Y_n)\}$.

With these definitions, we now summarize our main results, which are stated technically in Theorems 1–3. The first result (Theorem 1) establishes a general correspondence between the family of $f$-divergences and the family of optimized $\phi$-risks. In particular, let $R_\phi(Q)$ denote the *optimal $\phi$-risk*, meaning the $\phi$-risk obtained by optimizing over the discriminant $\gamma$ as follows:

$$R_\phi(Q) := \inf_{\gamma \in \Gamma} R_\phi(Q, \gamma).$$

In Theorem 1, we establish a precise correspondence between these optimal $\phi$-risks and the family of $f$-divergences. Theorem 1(a) addresses the forward direction of this correspondence (from $\phi$ to $f$); in particular, we show that any optimal $\phi$-risk can be written as $\mathbb{R}_\phi(Q) = -I_f(\mu, \pi)$, where $I_f$ is the divergence induced by a suitably chosen convex function $f$. We also specify a set of properties that any such function $f$ inherits from the surrogate loss $\phi$. Theorem 1(b) addresses the converse question: given an $f$-divergence, when can it be realized as an optimal $\phi$-risk? We provide a set of necessary and sufficient conditions on any such $f$-divergence and, moreover, specify a constructive procedure for determining *all* surrogate loss functions $\phi$ that induce the specified $f$-divergence.

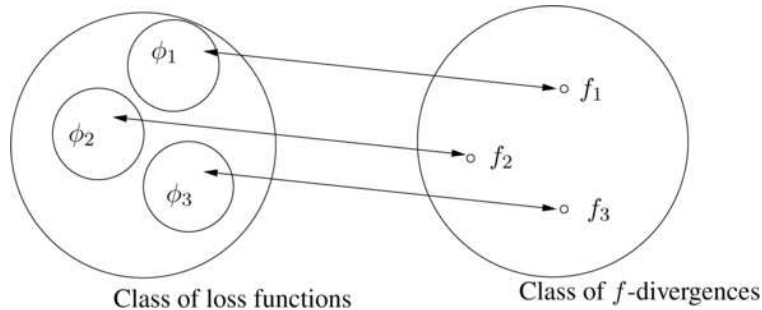

FIG. 1. *Illustration of the correspondence between f-divergences and loss functions. For each loss function $\phi$, there exists exactly one corresponding f-divergence such that the optimized $\phi$-risk is equal to the negative f-divergence. The reverse mapping is, in general, many-to-one.*

The relationship is illustrated in Figure 1; whereas each surrogate loss $\phi$ induces only one $f$-divergence, note that in general there are many surrogate loss functions that correspond to the same $f$-divergence. As particular examples of the general correspondence established in this paper, we show that the hinge loss corresponds to the variational distance, the exponential loss corresponds to the Hellinger distance, and the logistic loss corresponds to the capacitory discrimination distance.

This correspondence, in addition to its intrinsic interest as an extension of Blackwell's work, has a number of consequences. In Section 3, we show that it allows us to isolate a class of $\phi$-losses for which empirical risk minimization is consistent in the joint (quantizer and discriminant) estimation setting. Note in particular (e.g., from Blackwell's work) that the $f$-divergence associated with the 0–1 loss is the total variational distance. In Theorem 2, we specify a broader class of $\phi$-losses that induce the total variational distance and prove that, under standard technical conditions, an empirical risk minimization procedure based on any such $\phi$-risk is Bayes consistent. This broader class includes not only the nonconvex 0–1 loss, but also other convex and computationally tractable $\phi$-losses, including the hinge loss function that is well known in the context of support vector machines [6]. The key novelty in this result is that it applies to procedures that optimize simultaneously over the discriminant function $\gamma$ and the quantizer $Q$.

One interpretation of Theorem 2 is as specifying a set of surrogate loss functions $\phi$ that are universally equivalent to the 0–1 loss, in that empirical risk minimization procedures based on such $\phi$ yield classifier-quantizer pairs $(\gamma^*, Q^*)$ that achieve the Bayes risk. In Section 4, we explore this notion of universal equivalence between loss functions in more depth. In particular, we say that two loss functions $\phi_1$ and $\phi_2$ are *universally equivalent* if the optimal risks $R_{\phi_1}(Q)$ and $R_{\phi_2}(Q)$ induce the same ordering on quantizers, meaning the ordering $R_{\phi_1}(Q_a) \leq R_{\phi_1}(Q_b)$ holds if and only if $R_{\phi_2}(Q_a) \leq R_{\phi_2}(Q_b)$ for all quantizer pairs $Q_a$ and $Q_b$. Thus, the set of surrogate loss functions can be categorized into subclasses by this equivalence, where of particular interest are all surrogate loss functions that are equivalent (in the sense just defined) to the 0–1 loss. In Theorem 3, we provide an explicit and easily tested set of conditions for a $\phi$-risk to be equivalent to the 0–1 loss. One consequence is that procedures based on a $\phi$-risk outside of this family cannot be Bayes consistent for joint optimization of the discriminant $\gamma$ and quantizer $Q$. Thus, coupled with our earlier result in Theorem 2, we obtain a set of necessary and sufficient conditions on $\phi$-losses to be Bayes consistent in this joint estimation setting.

**2. Correspondence between $\phi$-loss and $f$-divergence.** Recall that in the setting of binary classification with $Q$ fixed, it is possible to give conditions on the class of surrogate loss functions (i.e., upper bounds on or approximations of the 0–1 loss) that yield Bayes consistency. In particular, Bartlett, Jordan and McAuliffe [2] have provided the following definition of a *classification-calibrated loss*.

DEFINITION 1. Define $\Phi_{a,b}(\alpha) = \phi(\alpha)a + \phi(-\alpha)b$. A loss function $\phi$ is *classification-calibrated* if for any $a, b \geq 0$ and $a \neq b$:

$$(3) \qquad \inf_{\{\alpha \in \mathbb{R} | \alpha(a-b) < 0\}} \Phi_{a,b}(\alpha) > \inf_{\{\alpha \in \mathbb{R} | \alpha(a-b) \geq 0\}} \Phi_{a,b}(\alpha).$$

The definition is essentially a pointwise form of a Fisher consistency condition that is appropriate for the binary classification setting. When $Q$ is fixed, this definition ensures that, under fairly general conditions, the decision rule $\gamma$ obtained by an empirical risk minimization procedure behaves equivalently to the Bayes optimal decision rule. Bartlett, Jordan and McAuliffe [2] also derived a simple lemma that characterizes classification-calibration for convex functions.

LEMMA 1. *Let $\phi$ be a convex function. Then $\phi$ is classification-calibrated if and only if it is differentiable at 0 and $\phi'(0) < 0$.*

For our purposes, we will find it useful to consider a somewhat more restricted definition of surrogate loss functions. In particular, we impose the following three conditions on any surrogate loss function $\phi : \mathbb{R} \to \mathbb{R} \cup \{+\infty\}$:

A1: $\phi$ is classification-calibrated;
A2: $\phi$ is continuous;
A3: Let $\alpha^* = \inf\{\alpha \in \mathbb{R} \cup \{+\infty\} | \phi(\alpha) = \inf \phi\}$. If $\alpha^* < +\infty$, then for any $\varepsilon > 0$,

$$\phi(\alpha^* - \varepsilon) \geq \phi(\alpha^* + \varepsilon). \tag{4}$$

The interpretation of assumption A3 is that one should penalize deviations away from $\alpha^*$ in the negative direction at least as strongly as deviations in the positive direction; this requirement is intuitively reasonable given the margin-based interpretation of $\alpha$. Moreover, this assumption is satisfied by all of the loss functions commonly considered in the literature; in particular, any decreasing function $\phi$ (e.g., hinge loss, logistic loss, exponential loss) satisfies this condition, as does the least squares loss (which is not decreasing). When $\phi$ is convex, assumption A1 is equivalent to requiring that $\phi$ be differentiable at 0 and $\phi'(0) < 0$. These facts also imply that the quantity $\alpha^*$ defined in assumption A3 is strictly positive. Finally, although $\phi$ is not defined for $-\infty$, we shall use the convention that $\phi(-\infty) = +\infty$.

In the following, we present the general relationship between optimal $\phi$-risks and $f$-divergences. The easier direction is to show that any $\phi$-risk induces a corresponding $f$-divergence. The $\phi$-risk can be written in the following way:

$$R_\phi(\gamma, Q) = \mathbb{E}\phi(Y\gamma(Z)) \tag{5a}$$

$$= \sum_z \phi(\gamma(z))\mu(z) + \phi(-\gamma(z))\pi(z). \tag{5b}$$

For a fixed mapping $Q$, the optimal $\phi$-risk has the form

$$R_\phi(Q) = \sum_{z \in \mathcal{Z}} \inf_\alpha (\phi(\alpha)\mu(z) + \phi(-\alpha)\pi(z))$$

$$= \sum_z \pi(z) \inf_\alpha \left( \phi(-\alpha) + \phi(\alpha) \frac{\mu(z)}{\pi(z)} \right).$$

For each $z$, define $u(z) := \frac{\mu(z)}{\pi(z)}$. With this notation, the function $\inf_\alpha(\phi(-\alpha) + \phi(\alpha)u)$ is concave as a function of $u$ (since the minimum of a collection of linear functions is concave). Thus, if we define

$$f(u) := -\inf_\alpha(\phi(-\alpha) + \phi(\alpha)u), \tag{6}$$

we obtain the relation

$$R_\phi(Q) = -I_f(\mu, \pi). \tag{7}$$

We have thus established the easy direction of the correspondence: given a loss function $\phi$, there exists an $f$-divergence for which the relation (7) holds. Furthermore, the convex function $f$ is given by the expression (6). Note that our argument does not require convexity of $\phi$.

We now consider the converse. Given a divergence $I_f(\mu, \pi)$ for some convex function $f$, does there exist a loss function $\phi$ for which $R_\phi(Q) = -I_f(\mu, \pi)$? In the theorem presented below, we answer this question in the affirmative. Moreover, we present a constructive result: we specify necessary and sufficient conditions under which there exist decreasing and convex surrogate loss functions for a given $f$-divergence, and we specify the form of all such loss functions.

Recall the notion of convex duality [18]: For a lower semicontinuous convex function $f : \mathbb{R} \to \mathbb{R} \cup \{\infty\}$, the conjugate dual $f^* : \mathbb{R} \to \mathbb{R} \cup \{\infty\}$ is defined as $f^*(u) = \sup_{v \in \mathbb{R}}(uv - f(v))$. Consider an intermediate function:

$$\Psi(\beta) = f^*(-\beta). \tag{8}$$

Define $\beta_1 := \inf\{\beta : \Psi(\beta) < +\infty\}$ and $\beta_2 := \inf\{\beta : \Psi(\beta) \leq \inf \Psi\}$. We are ready to state our first main result.

THEOREM 1. (a) *For any margin-based surrogate loss function $\phi$, there is an $f$-divergence such that $R_\phi(Q) = -I_f(\mu, \pi)$ for some lower semicontinuous convex function $f$.*

*In addition, if $\phi$ is a decreasing convex loss function that satisfies conditions* A1, A2 *and* A3, *then the following properties hold:*

  (i) $\Psi$ *is a decreasing and convex function;*
 (ii) $\Psi(\Psi(\beta)) = \beta$ *for all* $\beta \in (\beta_1, \beta_2)$;
(iii) *there exists a point* $u^* \in (\beta_1, \beta_2)$ *such that* $\Psi(u^*) = u^*$.

(b) *Conversely, if $f$ is a lower semicontinuous convex function satisfying all conditions (i)–(iii), there exists a decreasing convex surrogate loss $\phi$ that induces the $f$-divergence in the sense of equations (6) and (7).*

For proof of this theorem and additional properties, see Section 5.1.

REMARKS. (a) The existential statement in Theorem 1 can be strengthened to a constructive procedure, through which we specify how to obtain any $\phi$ loss function that induces a given $f$-divergence. Indeed, in the proof of Theorem 1(b) presented in Section 5.1, we prove that any decreasing surrogate loss function $\phi$ satisfying conditions A1–A3 that induces an $f$-divergence must be of the form

$$\phi(\alpha) = \begin{cases} u^*, & \text{if } \alpha = 0, \\ \Psi(g(\alpha + u^*)), & \text{if } \alpha > 0, \\ g(-\alpha + u^*), & \text{if } \alpha < 0, \end{cases} \tag{9}$$

where $g:[u^*, +\infty) \to \overline{\mathbb{R}}$ is some increasing continuous and convex function such that $g(u^*) = u^*$, and $g$ is right-differentiable at $u^*$ with $g'(u^*) > 0$.

(b) Another consequence of Theorem 1 is that any $f$-divergence can be obtained from a rather large set of surrogate loss functions; indeed, different such losses are obtained by varying the function $g$ in our constructive specification (9). In Section 2.1, we provide concrete examples of this constructive procedure and the resulting correspondences. For instance, we show that the variational distance corresponds to the 0–1 loss and the hinge loss, while the Hellinger distance corresponds to the exponential loss. Both divergences are also obtained from many less familiar loss functions.

(c) Although the correspondence has been formulated in the population setting, it is the basis of a constructive method for specifying a class of surrogate loss functions that yield a Bayes consistent estimation procedure. Indeed, in Section 3, we exploit this result to isolate a subclass of surrogate convex loss functions that yield Bayes-consistent procedures for joint $(\gamma, Q)$ minimization procedures. Interestingly, this class is a strict subset of the class of classification-calibrated loss functions, all of which yield Bayes-consistent estimation procedure in the standard classification setting (e.g., [2]). For instance, the class that we isolate contains the hinge loss, but *not* the exponential loss or the logistic loss functions. Finally, in Section 4, we show that, in a suitable sense, the specified subclass of surrogate loss functions is the only one that yields consistency for the joint $(\gamma, Q)$ estimation problem.

2.1. *Examples.* In this section, we describe various correspondences between $\phi$-losses and $f$-divergences that illustrate the claims of Theorem 1.

2.1.1. *0–1 loss, hinge loss and variational distance.* First, consider the 0–1 loss $\phi(\alpha) = \mathbb{I}[\alpha \leq 0]$. From equation (5b), the optimal discriminant function $\gamma$ takes the form $\gamma(z) = \operatorname{sign}(\mu(z) - \pi(z))$, so that the optimal Bayes risk is given by

$$R_{\text{Bayes}}(Q) = \sum_{z \in \mathcal{Z}} \min\{\mu(z), \pi(z)\}$$
$$= \tfrac{1}{2} - \tfrac{1}{2} \sum_{z \in \mathcal{Z}} |\mu(z) - \pi(z)| = \tfrac{1}{2}(1 - V(\mu, \pi)),$$

where $V(\mu, \pi)$ denotes the variational distance $V(\mu, \pi) := \sum_{z \in \mathcal{Z}} |\mu(z) - \pi(z)|$ between the two measures $\mu$ and $\pi$.

Now, consider the hinge loss function $\phi(\alpha) = \max\{0, 1 - \alpha\} = (1 - \alpha)_+$. In this case, a similar calculation yields $\gamma(z) = \operatorname{sign}(\mu(z) - \pi(z))$ as the optimal discriminant. The optimal risk for hinge loss thus takes the form:

$$R_{\text{hinge}}(Q) = \sum_{z \in \mathcal{Z}} 2 \min\{\mu(z), \pi(z)\} = 1 - \sum_{z \in \mathcal{Z}} |\mu(z) - \pi(z)| = 1 - V(\mu, \pi).$$

Thus, both the 0–1 loss and the hinge loss give rise to $f$-divergences of the form $f(u) = -c\min\{u,1\} + au + b$ for some constants $c > 0$ and $a, b$. Conversely, consider an $f$-divergence that is based on the function $f(u) = -2\min(u,1)$ for $u \geq 0$. Augmenting the definition by setting $f(u) = +\infty$ for $u < 0$, we use equation (9) to calculate $\Psi$:

$$\Psi(\beta) = f^*(-\beta) = \sup_{u \in \mathbb{R}}(-\beta u - f(u)) = \begin{cases} 0, & \text{if } \beta > 2, \\ 2 - \beta, & \text{if } 0 \leq \beta \leq 2, \\ +\infty, & \text{if } \beta < 0. \end{cases}$$

By inspection, we see that $u^* = 1$, where $u^*$ was defined in part (iii) of Theorem 1(a). If we set $g(u) = u$, then we recover the hinge loss $\phi(\alpha) = (1 - \alpha)_+$. On the other hand, choosing $g(u) = e^{u-1}$ leads to the loss

(10) $$\phi(\alpha) = \begin{cases} (2 - e^\alpha)_+, & \text{for } \alpha \leq 0, \\ e^{-\alpha}, & \text{for } \alpha > 0. \end{cases}$$

Note that the loss function obtained with this particular choice of $g$ is not convex, but our theory nonetheless guarantees that this non-convex loss still induces $f$ in the sense of equation (7). To ensure that $\phi$ is convex, we must choose $g$ to be an increasing convex function in $[1, +\infty)$ such that $g(u) = u$ for $u \in [1, 2]$. See Figure 2 for illustrations of some convex $\phi$ losses.

2.1.2. *Exponential loss and Hellinger distance.* Now, consider the exponential loss $\phi(\alpha) = \exp(-\alpha)$. In this case, a little calculation shows that the optimal discriminant is $\gamma(z) = \frac{1}{2}\log\frac{\mu(z)}{\pi(z)}$. The optimal risk for exponential loss is given by

$$R_{\exp}(Q) = \sum_{z \in \mathcal{Z}} 2\sqrt{\mu(z)\pi(z)} = 1 - \sum_{z \in \mathcal{Z}}(\sqrt{\mu(z)} - \sqrt{\pi(z)})^2 = 1 - 2h^2(\mu, \pi),$$

where $h(\mu, \pi) := \frac{1}{2}\sum_{z \in \mathcal{Z}}(\sqrt{\mu(z)} - \sqrt{\pi(z)})^2$ denotes the Hellinger distance between measures $\mu$ and $\pi$. Conversely, the Hellinger distance is equivalent to the negative of the Bhattacharyya distance, which is an $f$-divergence with $f(u) = -2\sqrt{u}$ for $u \geq 0$. Let us augment the definition of $f$ by setting $f(u) = +\infty$ for $u < 0$; doing so does not alter the Hellinger (or Bhattacharyya) distances. As before,

$$\Psi(\beta) = f^*(-\beta) = \sup_{u \in \mathbb{R}}(-\beta u - f(u)) = \begin{cases} 1/\beta, & \text{when } \beta > 0, \\ +\infty, & \text{otherwise.} \end{cases}$$

Thus, we see that $u^* = 1$. If we let $g(u) = u$, then a possible surrogate loss function that realizes the Hellinger distance takes the form:

$$\phi(\alpha) = \begin{cases} 1, & \text{if } \alpha = 0, \\ \dfrac{1}{\alpha + 1}, & \text{if } \alpha > 0, \\ -\alpha + 1, & \text{if } \alpha < 0. \end{cases}$$

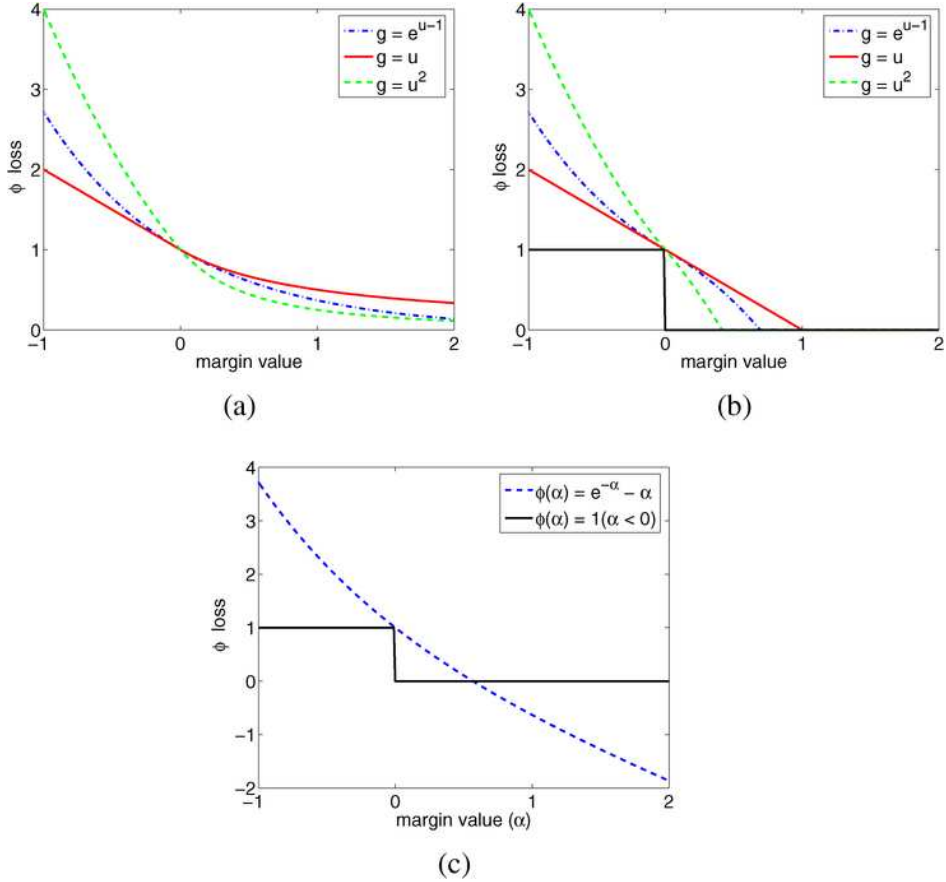

FIG. 2. *Panels* (a) *and* (b) *show examples of $\phi$ losses that induce the Hellinger distance and variational distance, respectively, based on different choices of the function $g$. Panel* (c) *shows a loss function that induces the symmetric KL divergence; for the purposes of comparison, the 0–1 loss is also plotted.*

On the other hand, if we set $g(u) = \exp(u-1)$, then we obtain the exponential loss $\phi(\alpha) = \exp(-\alpha)$. See Figure 2 for illustrations of these loss functions.

2.1.3. *Least squares loss and triangular discrimination distance.* Letting $\phi(\alpha) = (1-\alpha)^2$ be the least squares loss, the optimal discriminant is given by $\gamma(z) = \frac{\mu(z)-\pi(z)}{\mu(z)+\pi(z)}$. Thus, the optimal risk for least squares loss takes the form

$$R_{\mathrm{sqr}}(Q) = \sum_{z \in \mathcal{Z}} \frac{4\mu(z)\pi(z)}{\mu(z)+\pi(z)} = 1 - \sum_{z \in \mathcal{Z}} \frac{(\mu(z)-\pi(z))^2}{\mu(z)+\pi(z)} = 1 - \Delta(\mu, \pi),$$

where $\Delta(\mu, \pi)$ denotes the *triangular discrimination* distance [20]. Conversely, the triangular discriminatory distance is equivalent to the negative of the harmonic distance; it is an $f$-divergence with $f(u) = -\frac{4u}{u+1}$ for $u \geq 0$. Let us augment $f$ with $f(u) = +\infty$ for $u < 0$. We have

$$\Psi(\beta) = \sup_{u \in \mathbb{R}}(-\beta u - f(u)) = \begin{cases} (2 - \sqrt{\beta})^2, & \text{for } \beta \geq 0, \\ +\infty, & \text{otherwise.} \end{cases}$$

Clearly, $u^* = 1$. In this case, setting $g(u) = u^2$ gives the least square loss $\phi(\alpha) = (1 - \alpha)^2$.

2.1.4. *Logistic loss and capacitory discrimination distance.* Let $\phi(\alpha) = \log(1 + \exp(-\alpha))$ be the logistic loss. Then, $\gamma(z) = \log \frac{\mu(z)}{\pi(z)}$. As a result, the optimal risk for logistic loss is given by

$$R_{\log}(Q) = \sum_{z \in \mathcal{Z}} \mu(z) \log \frac{\mu(z) + \pi(z)}{\mu(z)} + \pi(z) \log \frac{\mu(z) + \pi(z)}{\pi(z)}$$

$$= \log 2 - KL\left(\mu \middle\| \frac{\mu + \pi}{2}\right) - KL\left(\pi \middle\| \frac{\mu + \pi}{2}\right) = \log 2 - C(\mu, \pi),$$

where $KL(U, V)$ denotes the Kullback–Leibler divergence between two measures $U$ and $V$, and $C(U, V)$ denotes the *capacitory discrimination* distance [20]. Conversely, the capacitory discrimination distance is equivalent to an $f$-divergence with $f(u) = -u \log \frac{u+1}{u} - \log(u+1)$, for $u \geq 0$. As before, augmenting this function with $f(u) = +\infty$ for $u < 0$, we have

$$\Psi(\beta) = \sup_{u \in \mathbb{R}}(-\beta u - f(u)) = \begin{cases} \beta - \log(e^\beta - 1), & \text{for } \beta \geq 0, \\ +\infty, & \text{otherwise.} \end{cases}$$

This representation shows that $u^* = \log 2$. If we choose $g(u) = \log(1 + \frac{e^u}{2})$, then we recover the logistic loss $\phi(\alpha) = \log[1 + \exp(-\alpha)]$.

2.1.5. *Another symmetrized Kullback–Leibler divergence.* Recall that both the KL divergences [i.e., $KL(\mu\|\pi)$ and $KL(\pi\|\mu)$] are asymmetric; therefore, Corollary 3 (see Section 5.1) implies that they are *not* realizable by any margin-based surrogate loss. However, a closely related functional is the *symmetric Kullback–Leibler* divergence [5]:

(11) $$KL_s(\mu, \pi) := KL(\mu\|\pi) + KL(\pi\|\mu).$$

It can be verified that this symmetrized KL divergence is an $f$-divergence, generated by the function $f(u) = -\log u + u \log u$ for $u \geq 0$, and $+\infty$ otherwise. Theorem 1 implies that it can be generated by surrogate loss functions of form (9), but the form of this loss function is not at all obvious. Therefore, in order to recover an explicit form for some $\phi$, we follow the constructive procedure outlined in the remarks following Theorem 1, first defining

$$\Psi(\beta) = \sup_{u \geq 0}\{-\beta u + \log u - u \log u\}.$$

In order to compute the value of this supremum, we take the derivative with respect to $u$ and set it to zero; doing so yields the zero-gradient condition $-\beta + 1/u - \log u - 1 = 0$. To capture this condition, we define a function $r:[0,+\infty) \to [-\infty,+\infty]$ via $r(u) = 1/u - \log u$. It is easy to see that $r$ is a strictly decreasing function whose range covers the whole real line; moreover, the zero-gradient condition is equivalent to $r(u) = \beta + 1$. We can thus write $\Psi(\beta) = u + \log u - 1$ where $u = r^{-1}(\beta + 1)$, or, equivalently,

$$\Psi(\beta) = r(1/u) - 1 = r\left(\frac{1}{r^{-1}(\beta+1)}\right) - 1.$$

It is straightforward to verify that the function $\Psi$ thus specified is strictly decreasing and convex with $\Psi(0) = 0$, and that $\Psi(\Psi(\beta)) = \beta$ for any $\beta \in \mathbb{R}$. Therefore, Theorem 1 allow us to specify the form of any convex surrogate loss function that generates the symmetric KL divergence; in particular, any such functions must be of the form (9):

$$\phi(\alpha) = \begin{cases} g(-\alpha), & \text{for } \alpha \leq 0, \\ \Psi(g(\alpha)), & \text{otherwise}, \end{cases}$$

where $g:[0,+\infty) \to [0,+\infty)$ is some increasing convex function satisfying $g(0) = 0$. As a particular example (and one that leads to a closed form expression for $\phi$), let us choose $g(u) = e^u + u - 1$. Doing so leads to the surrogate loss function

$$\phi(\alpha) = e^{-\alpha} - \alpha - 1,$$

as illustrated in Figure 2(c).

**3. Bayes consistency via surrogate losses.** As shown in Section 2.1.1, if we substitute the (nonconvex) 0–1 loss function into the linking equation (6), then we obtain the variational distance $V(\mu,\pi)$ as the $f$-divergence associated with the function $f(u) = \min\{u,1\}$. A bit more broadly, let us consider the subclass of $f$-divergences defined by functions of the form

(12) $$f(u) = -c\min\{u,1\} + au + b,$$

where $a, b$ and $c$ are scalars with $c > 0$. (For further examples of such losses, in addition to the 0–1 loss, see Section 2.1.) The main result of this section is that there exists a subset of surrogate losses $\phi$ associated with an $f$-divergence of the form (12) that, when used in the context of a risk minimization procedure for jointly optimizing $(\gamma, Q)$ pairs, yields a Bayes consistent method.

We begin by specifying some standard technical conditions under which our Bayes consistency result holds. Consider sequences of increasing compact function classes $\mathcal{C}_1 \subseteq \mathcal{C}_2 \subseteq \cdots \subseteq \Gamma$ and $\mathcal{D}_1 \subseteq \mathcal{D}_2 \subseteq \cdots \subseteq \mathcal{Q}$. Recall that $\Gamma$ denotes the class of all measurable functions from $\mathcal{Z} \to \mathbb{R}$, whereas $\mathcal{Q}$ is a constrained class of quantizer functions $Q$, with the restriction that $\mu$ and $\pi$ are strictly positive measures. Our analysis supposes that there exists an oracle that outputs an optimal solution to the minimization problem

$$(13) \quad \min_{(\gamma,Q)\in(\mathcal{C}_n,\mathcal{D}_n)} \hat{R}_\phi(\gamma, Q) = \min_{(\gamma,Q)\in(\mathcal{C}_n,\mathcal{D}_n)} \frac{1}{n} \sum_{i=1}^n \sum_{z \in \mathcal{Z}} \phi(Y_i \gamma(z)) Q(z|X_i),$$

and let $(\gamma_n^*, Q_n^*)$ denote one such solution. Let $R^*_{\text{Bayes}}$ denote the minimum Bayes risk achieved over the space of decision rules $(\gamma, Q) \in (\Gamma, \mathcal{Q})$:

$$(14) \quad R^*_{\text{Bayes}} := \inf_{(\gamma,Q)\in(\Gamma,\mathcal{Q})} R_{\text{Bayes}}(\gamma, Q).$$

We refer to the nonnegative quantity $R_{\text{Bayes}}(\gamma_n^*, Q_n^*) - R^*_{\text{Bayes}}$ as the *excess Bayes risk* of our estimation procedure. We say that such an estimation procedure is *universally consistent* if the excess Bayes risk converges to zero, that is, if under the (unknown) Borel probability measure $\mathbb{P}$ on $\mathcal{X} \times \mathcal{Y}$, we have

$$(15) \quad \lim_{n\to\infty} R_{\text{Bayes}}(\gamma_n^*, Q_n^*) = R^*_{\text{Bayes}} \qquad \text{in probability.}$$

In order to analyze the statistical behavior of this algorithm and to establish universal consistency for appropriate sequences $(\mathcal{C}_n, \mathcal{D}_n)$ of function classes, we follow a standard strategy of decomposing the Bayes error in terms of two types of errors:

- the *approximation error* associated with function classes $\mathcal{C}_n \subseteq \Gamma$, and $\mathcal{D}_n \subseteq \mathcal{Q}$:

$$(16) \quad \mathcal{E}_0(\mathcal{C}_n, \mathcal{D}_n) = \inf_{(\gamma,Q)\in(\mathcal{C}_n,\mathcal{D}_n)} \{R_\phi(\gamma, Q)\} - R^*_\phi,$$

  where $R^*_\phi := \inf_{(\gamma,Q)\in(\Gamma,\mathcal{Q})} R_\phi(\gamma, Q)$;
- the *estimation error* introduced by the finite sample size $n$:

$$(17) \quad \mathcal{E}_1(\mathcal{C}_n, \mathcal{D}_n) = \mathbb{E} \sup_{(\gamma,Q)\in(\mathcal{C}_n,\mathcal{D}_n)} |\hat{R}_\phi(\gamma, Q) - R_\phi(\gamma, Q)|,$$

  where the expectation is taken with respect to the (unknown) measure $\mathbb{P}^n(X, Y)$.

For asserting universal consistency, we impose the standard conditions:

(18) **Approximation condition:** $\quad \lim_{n\to\infty} \mathcal{E}_0(\mathcal{C}_n, \mathcal{D}_n) = 0.$

(19) **Estimation condition:** $\quad \lim_{n\to\infty} \mathcal{E}_1(\mathcal{C}_n, \mathcal{D}_n) = 0 \quad$ in probability.

**Conditions on loss function $\phi$:** Our consistency result applies to the class of surrogate losses that satisfy the following:

B1: $\phi$ is continuous, convex, and classification-calibrated;
B2: For each $n = 1, 2, \ldots$, we assume that

(20) $$M_n := \max_{y \in \{-1,+1\}} \sup_{(\gamma,Q) \in (\mathcal{C}_n, \mathcal{D}_n)} \sup_{z \in \mathcal{Z}} |\phi(y\gamma(z))| < +\infty.$$

With this set-up, the following theorem ties together the Bayes error with the approximation error and estimation error and provides sufficient conditions for universal consistency for a *suitable subclass* of surrogate loss functions.

THEOREM 2. *Consider an estimation procedure of the form (13), using a surrogate loss $\phi$. Recall the prior probabilities $p = \mathbb{P}(Y=1)$ and $q = \mathbb{P}(Y=-1)$. For any surrogate loss $\phi$ satisfying conditions* B1 *and* B2 *and inducing an $f$-divergence of the form (12) for any $c > 0$, and for $a, b$ such that $(a-b)(p-q) \geq 0$, we have:*

(a) *For any Borel probability measure $\mathbb{P}$, there holds, with probability at least $1 - \delta$:*

$$R_{\text{Bayes}}(\gamma_n^*, Q_n^*) - R_{\text{Bayes}}^*$$
$$\leq \frac{2}{c} \left\{ 2\mathcal{E}_1(\mathcal{C}_n, \mathcal{D}_n) + \mathcal{E}_0(\mathcal{C}_n, \mathcal{D}_n) + 2M_n \sqrt{2\frac{\ln(2/\delta)}{n}} \right\}.$$

(b) *Universal Consistency: For function classes satisfying the approximation (18) and estimation conditions (19), the estimation procedure (13) is universally consistent:*

(21) $\quad \lim_{n\to\infty} R_{\text{Bayes}}(\gamma_n^*, Q_n^*) = R_{\text{Bayes}}^* \quad$ *in probability.*

REMARKS. (i) Note that both the approximation and the estimation errors are with respect to the $\phi$-loss, but the theorem statement refers to the excess Bayes risk. Since the analysis of approximation and estimation conditions such as those in equation (18) and (19) is a standard topic in statistical learning, we will not discuss it further here. We note that our previous work analyzed the estimation error for certain kernel classes [15].

(ii) It is worth pointing out that in order for our result to be applicable to an *arbitrary* constrained class of $\mathcal{Q}$ for which $\mu$ and $\pi$ are strictly positive measures, we need the additional constraint that $(a-b)(p-q) \geq 0$, where $a, b$ are scalars in the $f$-divergence (12) and $p, q$ are the unknown prior probabilities. Intuitively, this requirement is needed to ensure that the approximation error due to varying $Q$ within $\mathcal{Q}$ dominates the approximation error due to varying $\gamma$ (because the optimal $\gamma$ is determined only after $Q$) for arbitrary $\mathcal{Q}$. Since $p$ and $q$ are generally unknown, the only $f$-divergences that are practically useful are the ones for which $a = b$. One such $\phi$ is the hinge loss, which underlies the support vector machine.

Finally, we note that the proof of Theorem 2 relies on an auxiliary result that is of independent interest. In particular, we prove that for any function classes $\mathcal{C}$ and $\mathcal{D}$, for certain choice of surrogate loss $\phi$, the excess $\phi$-risk is related to the excess Bayes risk as follows.

LEMMA 2. *Let $\phi$ be a surrogate loss function satisfying all conditions specified in Theorem 2. Then, for any classifier-quantizer pair $(\gamma, Q) \in (\mathcal{C}, \mathcal{D})$, we have*

$$(22) \qquad \frac{c}{2}[R_{\text{Bayes}}(\gamma, Q) - R^*_{\text{Bayes}}] \leq R_\phi(\gamma, Q) - R^*_\phi.$$

This result (22) demonstrates that in order to achieve joint Bayes consistency—that is, in order to drive the excess Bayes risk to zero, while optimizing over the pair $(\gamma, Q)$—it suffices to drive the excess $\phi$-risk to zero.

**4. Comparison between loss functions.** We have studied a broad class of loss functions corresponding to $f$-divergences of the form (12) in Theorem 1. A subset of this class in turn yields Bayes consistency for the estimation procedure (13) as shown in Theorem 2. A natural question is, are there any other surrogate loss functions that also yield Bayes consistency?

A necessary condition for achieving Bayes consistency using estimation procedure (13) is that the constrained minimization over surrogate $\phi$-risks should yield a $(Q, \gamma)$ pair that minimizes the expected 0–1 loss subject to the same constraints. In this section, we show that only surrogate loss functions that induce $f$-divergence of the form (12) can actually satisfy this property. We establish this result by developing a general way of comparing different loss functions. In particular, by exploiting the correspondence between surrogate losses and $f$-divergences, we are able to compare surrogate losses in terms of their corresponding $f$-divergences.

4.1. *Connection between 0–1 loss and $f$-divergences.* The connection between $f$-divergences and 0–1 loss that we develop has its origins in seminal work on comparison of experiments by Blackwell and others [3, 4, 5]. In particular, we give the following definition.

DEFINITION 2. The quantizer $Q_1$ dominates $Q_2$ if $R_{\text{Bayes}}(Q_1) \leq R_{\text{Bayes}}(Q_2)$ for any choice of prior probability $q = \mathbb{P}(Y = -1) \in (0,1)$.

Recall that a choice of quantizer design $Q$ induces two conditional distributions, say $P(Z|Y = 1) \sim P_1$ and $P(Z|Y = -1) \sim P_{-1}$. From here onward, we use $P_{-1}^Q$ and $P_1^Q$ to denote the fact that both $P_{-1}$ and $P_1$ are determined by the specific choice of $Q$. By "parameterizing" the decision-theoretic criterion in terms of loss function $\phi$ and establishing a precise correspondence between $\phi$ and the $f$-divergence, we obtain an arguably simpler proof of the classical theorem [3, 4] that relates 0–1 loss to $f$-divergences.

PROPOSITION 1 [3, 4]. *For any two quantizer designs $Q_1$ and $Q_2$, the following statements are equivalent:*

(a) $Q_1$ *dominates* $Q_2$ *[i.e., $R_{\text{Bayes}}(Q_1) \leq R_{\text{Bayes}}(Q_2)$ for any prior probability $q \in (0,1)$];*
(b) $I_f(P_1^{Q_1}, P_{-1}^{Q_1}) \geq I_f(P_1^{Q_2}, P_{-1}^{Q_2})$, *for all functions $f$ of the form $f(u) = -\min(u,c)$ for some $c > 0$;*
(c) $I_f(P_1^{Q_1}, P_{-1}^{Q_1}) \geq I_f(P_1^{Q_2}, P_{-1}^{Q_2})$, *for all convex functions $f$.*

PROOF. We first establish the equivalence (a) $\Leftrightarrow$ (b). By the correspondence between 0–1 loss and an $f$-divergence with $f(u) = -\min(u,1)$, we have $R_{\text{Bayes}}(Q) = -I_f(\mu, \pi) = -I_{f_q}(P_1, P_{-1})$, where $f_q(u) := qf(\frac{1-q}{q}u) = -(1-q)\min(u, \frac{q}{1-q})$. Hence, (a) $\Leftrightarrow$ (b).

Next, we prove the equivalence (b) $\Leftrightarrow$ (c). The implication (c) $\Rightarrow$ (b) is immediate. Considering the reverse implication (b) $\Rightarrow$ (c), we note that any convex function $f(u)$ can be uniformly approximated over a bounded interval as a sum of a linear function and $-\sum_k \alpha_k \min(u, c_k)$, where $\alpha_k > 0, c_k > 0$ for all $k$. For a linear function $f$, $I_f(P_{-1}, P_1)$ does not depend on $P_{-1}, P_1$. Using these facts, (c) follows easily from (b). □

COROLLARY 1. *The quantizer $Q_1$ dominates $Q_2$ if and only if $R_\phi(Q_1) \leq R_\phi(Q_2)$ for any loss function $\phi$.*

PROOF. By Theorem 1(a), we have $R_\phi(Q) = -I_f(\mu, \pi) = -I_{f_q}(P_1, P_{-1})$, from which the corollary follows, using Proposition 1. □

Corollary 1 implies that if $R_\phi(Q_1) \leq R_\phi(Q_2)$ for some loss function $\phi$, then $R_{\text{Bayes}}(Q_1) \leq R_{\text{Bayes}}(Q_2)$ for some set of prior probabilities on the hypothesis space. This implication justifies the use of a given surrogate loss function $\phi$ in place of the 0–1 loss for *some* prior probability; however, for a given prior probability, it gives no guidance on how to choose $\phi$. Moreover, the prior probabilities on the label $Y$ are typically unknown in many applications. In such a setting, Blackwell's notion of $Q_1$ dominating $Q_2$ has limited usefulness. With this motivation in mind, the following section is devoted to development of a more stringent method for assessing equivalence between loss functions.

4.2. *Universal equivalence.* Suppose that the loss functions $\phi_1$ and $\phi_2$ realize the $f$-divergences associated with the convex functions $f_1$ and $f_2$, respectively. We then have the following definition.

DEFINITION 3. The surrogate loss functions $\phi_1$ and $\phi_2$ are *universally equivalent*, denoted by $\phi_1 \stackrel{u}{\approx} \phi_2$, if for any $\mathbb{P}(X,Y)$ and quantization rules $Q_1, Q_2$, there holds:
$$R_{\phi_1}(Q_1) \leq R_{\phi_1}(Q_2) \Leftrightarrow R_{\phi_2}(Q_1) \leq R_{\phi_2}(Q_2).$$
In terms of the corresponding $f$-divergences, this relation is denoted by $f_1 \stackrel{u}{\approx} f_2$.

Observe that this definition is very stringent, in that it requires that the ordering between optimal $\phi_1$ and $\phi_2$ risks holds for all probability distributions $\mathbb{P}$ on $\mathcal{X} \times \mathcal{Y}$. However, this stronger notion of equivalence is needed for nonparametric approaches to classification, in which the underlying distribution $\mathbb{P}$ is only weakly constrained.

The following result provides necessary and sufficient conditions for two $f$-divergences to be universally equivalent.

THEOREM 3. *Let $f_1$ and $f_2$ be continuous, nonlinear and convex functions on $[0, +\infty) \to \mathbb{R}$. Then, $f_1 \stackrel{u}{\approx} f_2$ if and only if $f_1(u) = cf_2(u) + au + b$ for some constants $c > 0$ and $a, b$.*

An important special case is when one of the $f$-divergences is the variational distance. In this case, we have the following.

COROLLARY 2. (a) *All $f$-divergences based on continuous convex $f : [0, +\infty) \to \infty$ that are universally equivalent to the variational distance have the form*

(23)	$$f(u) = -c\min(u, 1) + au + b \quad \textit{for some } c > 0.$$

(b) *The 0–1 loss is universally equivalent only to those loss functions whose corresponding $f$-divergence is based on a function of the form (23).*

The above result establishes that only those surrogate loss functions corresponding to the variational distance yield universal consistency in a strong sense, meaning for any underlying $\mathbb{P}$ and a constrained class of quantization rules.

**5. Proofs.** In this section, we provide detailed proofs of our main results, as well as some auxiliary results.

5.1. *Proofs of Theorem 1 and auxiliary properties.* Our proof proceeds via connecting some intermediate functions. First, let us define, for each $\beta$, the inverse mapping

$$\phi^{-1}(\beta) := \inf\{\alpha : \phi(\alpha) \leq \beta\}, \tag{24}$$

where $\inf \varnothing := +\infty$. The following result summarizes some useful properties of $\phi^{-1}$.

LEMMA 3. *Suppose that $\phi$ is a convex loss satisfying assumptions* A1, A2 *and* A3.

(a) *For all $\beta \in \mathbb{R}$ such that $\phi^{-1}(\beta) < +\infty$, the inequality $\phi(\phi^{-1}(\beta)) \leq \beta$ holds. Furthermore, equality occurs when $\phi$ is continuous at $\phi^{-1}(\beta)$.*
(b) *The function $\phi^{-1} : \mathbb{R} \to \overline{\mathbb{R}}$ is strictly decreasing and convex.*

Using the function $\phi^{-1}$, we define a new function $\tilde{\Psi} : \mathbb{R} \to \overline{\mathbb{R}}$ by

$$\tilde{\Psi}(\beta) := \begin{cases} \phi(-\phi^{-1}(\beta)), & \text{if } \phi^{-1}(\beta) \in \mathbb{R}, \\ +\infty, & \text{otherwise.} \end{cases} \tag{25}$$

Note that the domain of $\tilde{\Psi}$ is $\text{Dom}(\tilde{\Psi}) = \{\beta \in \mathbb{R} : \phi^{-1}(\beta) \in \mathbb{R}\}$. Now, define

$$(26) \quad \tilde{\beta}_1 := \inf\{\beta : \tilde{\Psi}(\beta) < +\infty\} \quad \text{and} \quad \tilde{\beta}_2 := \inf\{\beta : \tilde{\Psi}(\beta) = \inf \tilde{\Psi}\}.$$

It is simple to check that $\inf \phi = \inf \tilde{\Psi} = \phi(\alpha^*)$, and $\tilde{\beta}_1 = \phi(\alpha^*)$, $\tilde{\beta}_2 = \phi(-\alpha^*)$. Furthermore, by construction, we have $\tilde{\Psi}(\tilde{\beta}_2) = \phi(\alpha^*) = \tilde{\beta}_1$, as well as $\tilde{\Psi}(\tilde{\beta}_1) = \phi(-\alpha^*) = \tilde{\beta}_2$. The following properties of $\tilde{\Psi}$ are particularly useful for our main results.

LEMMA 4. *Suppose that $\phi$ is a convex loss satisfying assumptions* A1, A2 *and* A3. *We have:*

(a) *$\tilde{\Psi}$ is strictly decreasing in the interval $(\tilde{\beta}_1, \tilde{\beta}_2)$. If $\phi$ is decreasing, then $\tilde{\Psi}$ is also decreasing in $(-\infty, +\infty)$. In addition, $\tilde{\Psi}(\beta) = +\infty$ for $\beta < \tilde{\beta}_1$.*

(b) $\tilde{\Psi}$ *is convex in* $(-\infty, \tilde{\beta}_2]$. *If* $\phi$ *is a decreasing function, then* $\tilde{\Psi}$ *is convex in* $(-\infty, +\infty)$.

(c) $\tilde{\Psi}$ *is lower semi-continuous, and continuous in its domain.*

(d) *For any* $\alpha \geq 0$, $\phi(\alpha) = \tilde{\Psi}(\phi(-\alpha))$. *In particular, there exists* $u^* \in (\tilde{\beta}_1, \tilde{\beta}_2)$ *such that* $\tilde{\Psi}(u^*) = u^*$.

(e) *The function* $\tilde{\Psi}$ *satisfies* $\tilde{\Psi}(\tilde{\Psi}(\beta)) \leq \beta$ *for all* $\beta \in Dom(\tilde{\Psi})$. *Moreover, if* $\phi$ *is a continuous function on its domain* $\{\alpha \in \mathbb{R} | \phi(\alpha) < +\infty\}$, *then* $\tilde{\Psi}(\tilde{\Psi}(\beta)) = \beta$ *for all* $\beta \in (\tilde{\beta}_1, \tilde{\beta}_2)$.

Let us proceed to part (a) of the theorem. The statement for general $\phi$ has already proved in the derivation preceding the theorem statement. Now, supposing that a decreasing convex surrogate loss $\phi$ satisfies assumptions A1, A2 and A3, then

$$f(u) = -\inf_{\alpha \in \mathbb{R}} (\phi(-\alpha) + \phi(\alpha)u)$$
$$= -\inf_{\{\alpha, \beta | \phi^{-1}(\beta) \in \mathbb{R}, \phi(\alpha) = \beta\}} (\phi(-\alpha) + \beta u).$$

For $\beta$ such that $\phi^{-1}(\beta) \in \mathbb{R}$, there might be more than one $\alpha$ such that $\phi(\alpha) = \beta$. However, our assumption (4) ensures that $\alpha = \phi^{-1}(\beta)$ results in minimum $\phi(-\alpha)$. Hence,

$$f(u) = -\inf_{\beta : \phi^{-1}(\beta) \in \mathbb{R}} (\phi(-\phi^{-1}(\beta)) + \beta u) = -\inf_{\beta \in \mathbb{R}} (\beta u + \tilde{\Psi}(\beta))$$
$$= \sup_{\beta \in \mathbb{R}} (-\beta u - \tilde{\Psi}(\beta)) = \tilde{\Psi}^*(-u).$$

By Lemma 4(b), the fact that $\phi$ is decreasing implies that $\tilde{\Psi}$ is convex. By convex duality and the lower semicontinuity of $\tilde{\Psi}$ (from Lemma 4(c)), we can also write

(27) $$\tilde{\Psi}(\beta) = \tilde{\Psi}^{**}(\beta) = f^*(-\beta).$$

Thus, $\tilde{\Psi}$ is identical to the function $\Psi$ defined in equation (8). The proof of part (a) is complete, thanks to Lemma 4. Furthermore, it can be shown that $\phi$ must have the form (9). Indeed, from Lemma 4(d), we have $\Psi(\phi(0)) = \phi(0) \in (\beta_1, \beta_2)$. As a consequence, $u^* := \phi(0)$ satisfies the relation $\Psi(u^*) = u^*$. Since $\phi$ is decreasing and convex on the interval $(-\infty, 0]$, for any $\alpha \geq 0$, we can write

$$\phi(-\alpha) = g(\alpha + u^*),$$

where $g$ is some increasing continuous and convex function. From Lemma 4(d), we have $\phi(\alpha) = \Psi(\phi(-\alpha)) = \Psi(g(\alpha + u^*)$ for $\alpha \geq 0$. To ensure the continuity at 0, there holds $u^* = \phi(0) = g(u^*)$. To ensure that $\phi$ is classification-calibrated, we require that $\phi$ be differentiable at 0 and $\phi'(0) < 0$. These conditions in turn imply that $g$ must be right-differentiable at $u^*$, with $g'(u^*) > 0$.

Let us turn to part (b) of the theorem. Since $f$ is lower semicontinuous by assumption, convex duality allows us to write
$$\begin{aligned} f(u) &= f^{**}(u) = \Psi^*(-u) \\ &= \sup_{\beta \in \mathbb{R}}(-\beta u - \Psi(\beta)) = -\inf_{\beta \in \mathbb{R}}(\beta u + \Psi(\beta)). \end{aligned}$$

Note that $\Psi$ is lower semicontinuous and convex by definition. To prove that any surrogate loss $\phi$ of form (9) (along with conditions A1–A3) must induce $f$-divergences in the sense of equation (6) [and thus equation (7)], it remains to show that $\phi$ is linked to $\Psi$ via the relation

(28) $$\tilde{\Psi} \equiv \Psi.$$

Since $\Psi$ is assumed to be a decreasing function, the function $\phi$ defined in (9) is also a decreasing function. Using the fixed point $u^* \in (\beta_1, \beta_2)$ of function $\Psi$, we divide our analysis into three cases:

- For $\beta \geq u^*$, there exists $\alpha \geq 0$ such that $g(\alpha + u^*) = \beta$. Choose the largest such $\alpha$. From our definition of $\phi$, $\phi(-\alpha) = \beta$. Thus, $\phi^{-1}(\beta) = -\alpha$. It follows that $\tilde{\Psi}(\beta) = \phi(-\phi^{-1}(\beta)) = \phi(\alpha) = \Psi(g(\alpha + u^*)) = \Psi(\beta)$.
- For $\beta < \beta_1$, then $\Psi(\beta) = +\infty$. It can also be verified that $\tilde{\Psi}(\beta) = +\infty$.
- Lastly, for $\beta_1 \leq \beta < u^* < \beta_2$, there exists $\alpha > 0$ such that $g(\alpha + u^*) \in (u^*, \beta_2)$ and $\beta = \Psi(g(\alpha + u^*))$, which implies that $\beta = \phi(\alpha)$ from our definition. Choose the smallest $\alpha$ that satisfies these conditions. Then, $\phi^{-1}(\beta) = \alpha$, and it follows that $\tilde{\Psi}(\beta) = \phi(-\phi^{-1}(\beta)) = \phi(-\alpha) = g(\alpha + u^*) = \Psi(\Psi(g(\alpha + u^*))) = \Psi(\beta)$, where we have used the fact that $g(\alpha + u^*) \in (\beta_1, \beta_2)$.

The proof of Theorem 1 is complete.

5.1.1. *Some additional properties.* In the remainder of this section we present several useful properties of surrogate losses and $f$-divergences. Although Theorem 1 provides one set of conditions for an $f$-divergence to be realized by some surrogate loss $\phi$, as well as a constructive procedure for finding all such loss functions, the following result provides a related set of conditions that can be easier to verify. We say that an $f$-divergence is *symmetric* if $I_f(\mu, \pi) = I_f(\pi, \mu)$ for any measures $\mu$ and $\pi$. With this definition, we have the following.

COROLLARY 3.  *Suppose that $f : [0, +\infty) \to \mathbb{R}$ is a continuous and convex function. The following are equivalent:*

(a) *The function $f$ is realizable by some surrogate loss function $\phi$ (via Theorem 1).*

(b) *The $f$-divergence $I_f$ is symmetric.*
(c) *For any $u > 0$, $f(u) = uf(1/u)$.*

PROOF. (a) $\Rightarrow$ (b): From Theorem 1(a), we have the representation $R_\phi(Q) = -I_f(\mu, \pi)$. Alternatively, we can write

$$R_\phi(Q) = \sum_z \mu(z) \min_\alpha \left( \phi(\alpha) + \phi(-\alpha) \frac{\pi(z)}{\mu(z)} \right)$$
$$= -\sum_z \mu(z) f\left( \frac{\pi(z)}{\mu(z)} \right),$$

which is equal to $-I_f(\pi, \mu)$, thereby showing that the $f$-divergence is symmetric.

(b) $\Rightarrow$ (c): By assumption, the following relation holds for any measures $\mu$ and $\pi$:

$$\sum_z \pi(z) f(\mu(z)/\pi(z)) = \sum_z \mu(z) f(\pi(z)/\mu(z)). \tag{29}$$

Take any instance of $z = l \in \mathcal{Z}$, and consider measures $\mu'$ and $\pi'$, which are defined on the space $\mathcal{Z} - \{l\}$ such that $\mu'(z) = \mu(z)$ and $\pi'(z) = \pi(z)$ for all $z \in \mathcal{Z} - \{l\}$. Since condition (29) also holds for $\mu'$ and $\pi'$, it follows that

$$\pi(z) f(\mu(z)/\pi(z)) = \mu(z) f(\pi(z)/\mu(z))$$

for all $z \in \mathcal{Z}$ and any $\mu$ and $\pi$. Hence, $f(u) = uf(1/u)$ for any $u > 0$.

(c) $\Rightarrow$ (a): It suffices to show that all sufficient conditions specified by Theorem 1 are satisfied.

Since any $f$-divergence is defined by applying $f$ to a likelihood ratio [see definition (2)], we can assume $f(u) = +\infty$ for $u < 0$ without loss of generality. Since $f(u) = uf(1/u)$ for any $u > 0$, it can be verified using subdifferential calculus [8] that for any $u > 0$, there holds

$$\partial f(u) = f(1/u) + \partial f(1/u) \frac{-1}{u}. \tag{30}$$

Given some $u > 0$, consider any $v_1 \in \partial f(u)$. Combined with equation (30) and the equality $f(u) = uf(1/u)$, we have

$$f(u) - v_1 u \in \partial f(1/u). \tag{31}$$

By definition of conjugate duality, $f^*(v_1) = v_1 u - f(u)$.

Letting $\Psi(\beta) = f^*(-\beta)$ as in Theorem 1, we have

$$\Psi(\Psi(-v_1)) = \Psi(f^*(v_1)) = \Psi(v_1 u - f(u))$$
$$= f^*(f(u) - v_1 u) = \sup_{\beta \in \mathbb{R}} (\beta f(u) - \beta v_1 u - f(\beta)).$$

Note that from equation (31), the supremum is achieved at $\beta = 1/u$, so that we have $\Psi(\Psi(-v_1)) = -v_1$ for any $v_1 \in \partial f(u)$ for $u > 0$. In other words, $\Psi(\Psi(\beta)) = \beta$ for any $\beta \in \{-\partial f(u), u > 0\}$. Convex duality and the definition $\Psi(\beta) = f^*(-\beta)$ imply that $\beta \in -\partial f(u)$ for some $u > 0$ if and only if $-u \in \partial \Psi(\beta)$ for some $u > 0$. This condition on $\beta$ is equivalent to the subdifferential $\partial \Psi(\beta)$ containing some negative value, which is satisfied by any $\beta \in (\beta_1, \beta_2)$, so that $\Psi(\Psi(\beta)) = \beta$ for $\beta \in (\beta_1, \beta_2)$. In addition, since $f(u) = +\infty$ for $u < 0$, $\Psi$ is a decreasing function. Now, as an application of Theorem 1, we conclude that $I_f$ is realizable by some (decreasing) surrogate loss function. $\square$

The following result establishes a link between (un)boundedness and the properties of the associated $f$.

COROLLARY 4. *Assume that $\phi$ is a decreasing (continuous convex) loss function corresponding to an $f$-divergence, where $f$ is a continuous convex function that is bounded from below by an affine function. Then, $\phi$ is unbounded from below if and only if $f$ is 1-coercive, that is, $f(x)/\|x\| \to +\infty$ as $\|x\| \to \infty$.*

PROOF. $\phi$ is unbounded from below if and only if $\Psi(\beta) = \phi(-\phi^{-1}(\beta)) \in \mathbb{R}$ for all $\beta \in \mathbb{R}$, which is equivalent to the dual function $f(\beta) = \Psi^*(-\beta)$ being 1-coercive cf. [8]. $\square$

Consequentially, for any decreasing and lower-bounded $\phi$ loss (which includes the hinge, logistic and exponential losses), the associated $f$-divergence is *not* 1-coercive. Other interesting $f$-divergences such as the *symmetric* KL divergence considered in [5] are 1-coercive, meaning that any associated surrogate loss $\phi$ cannot be bounded below.

5.2. *Proof of Theorem 2.* First let us prove Lemma 2:

PROOF. Since $\phi$ has form (9), it is easy to check that $\phi(0) = (c-a-b)/2$. Now, note that

$$\begin{aligned}
R_{\text{Bayes}}(\gamma, Q) - R_{\text{Bayes}}^* &= R_{\text{Bayes}}(\gamma, Q) - R_{\text{Bayes}}(Q) + R_{\text{Bayes}}(Q) - R_{\text{Bayes}}^* \\
&= \sum_{z \in \mathcal{Z}} \pi(z) \mathbb{I}(\gamma(z) > 0) + \mu(z) \mathbb{I}(\gamma(z) < 0) \\
&\quad - \min\{\mu(z), \pi(z)\} + R_{\text{Bayes}}(Q) - R_{\text{Bayes}}^* \\
&= \sum_{z:(\mu(z)-\pi(z))\gamma(z)<0} |\mu(z) - \pi(z)| + R_{\text{Bayes}}(Q) - R_{\text{Bayes}}^*.
\end{aligned}$$

In addition,
$$R_\phi(\gamma, Q) - R_\phi^* = R_\phi(\gamma, Q) - R_\phi(Q) + R_\phi(Q) - R_\phi^*.$$

By Theorem 1(a),

$$R_\phi(Q) - R_\phi^* = -I_f(\mu,\pi) - \inf_{Q\in\mathcal{Q}}(-I_f(\mu,\pi))$$

$$= c\sum_{z\in\mathcal{Z}}\min\{\mu(z),\pi(z)\} - \inf_{Q\in\mathcal{Q}} c\sum_{z\in\mathcal{Z}}\min\{\mu(z),\pi(z)\}$$

$$= c(R_{\text{Bayes}}(Q) - R_{\text{Bayes}}^*).$$

Therefore, the lemma will be immediate once we can show that

$$\frac{c}{2}\sum_{z:(\mu(z)-\pi(z))\gamma(z)<0}|\mu(z)-\pi(z)| \leq R_\phi(\gamma,Q) - R_\phi(Q)$$

(32)
$$= \sum_{z\in\mathcal{Z}}\pi(z)\phi(-\gamma(z)) + \mu(z)\phi(\gamma(z))$$

$$- c\min\{\mu(z),\pi(z)\} + ap + bq.$$

It is easy to check that for any $z \in \mathcal{Z}$ such that $(\mu(z) - \pi(z))\gamma(z) < 0$, there holds

(33) $$\pi(z)\phi(-\gamma(z)) + \mu(z)\phi(\gamma(z)) \geq \pi(z)\phi(0) + \mu(z)\phi(0).$$

Indeed, without loss of generality, suppose $\mu(z) > \pi(z)$. Since $\phi$ is classification-calibrated, the convex function (with respect to $\alpha$) $\pi(z)\phi(-\alpha) + \mu(z)\phi(\alpha)$ achieves its minimum at some $\alpha \geq 0$. Hence, for any $\alpha \leq 0$, $\pi(z)\phi(-\alpha) + \mu(z)\phi(\alpha) \geq \pi(z)\phi(0) + \mu(z)\phi(0)$. Hence, the statement (33) is proven. The RHS of equation (32) is lower bounded by

$$\sum_{z:(\mu(z)-\pi(z))\gamma(z)<0}(\pi(z)+\mu(z))\phi(0) - c\min\{\mu(z),\pi(z)\} + ap + bq$$

$$= \sum_{z:(\mu(z)-\pi(z))\gamma(z)<0}(\pi(z)+\mu(z))\frac{c-a-b}{2} - c\min\{\mu(z),\pi(z)\}$$

$$+ ap + bq$$

$$\geq \frac{c}{2}\sum_{z:(\mu(z)-\pi(z))\gamma(z)<0}|\mu(z)-\pi(z)| - (a+b)(p+q)/2 + ap + bq$$

$$= \frac{c}{2}\sum_{z:(\mu(z)-\pi(z))\gamma(z)<0}|\mu(z)-\pi(z)| + \frac{1}{2}(a-b)(p-q)$$

$$\geq \frac{c}{2}\sum_{z:(\mu(z)-\pi(z))\gamma(z)<0}|\mu(z)-\pi(z)|.$$

This completes the proof of the lemma. □

We are now equipped to prove Theorem 2. For part (a), first observe that the value of $\sup_{\gamma \in \mathcal{C}_n, Q \in \mathcal{D}_n} |\hat{R}_\phi(\gamma, Q) - R_\phi(\gamma, Q)|$ varies by at most $2M_n/n$ if one changes the values of $(X_i, Y_i)$ for some index $i \in \{1, \ldots, n\}$. Hence, applying McDiarmid's inequality yields concentration around the expected value [14], or (alternatively stated) we have that, with probability at least $1 - \delta$,

$$(34) \quad \left| \sup_{\gamma \in \mathcal{C}_n, Q \in \mathcal{D}_n} |\hat{R}_\phi(\gamma, Q) - R_\phi(\gamma, Q)| - \mathcal{E}_1(\mathcal{C}_n, \mathcal{D}_n) \right| \leq M_n \sqrt{2\ln(1/\delta)/n}.$$

Suppose that $R_\phi(\gamma, Q)$ attains its minimum over the compact subset $(\mathcal{C}_n, \mathcal{D}_n)$ at $(\gamma_n^\dagger, Q_n^\dagger)$. Then, using Lemma 2, we have

$$\frac{c}{2}(R_{\text{Bayes}}(\gamma_n^*, Q_n^*) - R_{\text{Bayes}}^*) \leq R_\phi(\gamma_n^*, Q_n^*) - R_\phi^*$$
$$= R_\phi(\gamma_n^*, Q_n^*) - R_\phi(\gamma_n^\dagger, Q_n^\dagger) + R_\phi(\gamma_n^\dagger, Q_n^\dagger) - R_\phi^*$$
$$= R_\phi(\gamma_n^*, Q_n^*) - R_\phi(\gamma_n^\dagger, Q_n^\dagger) + \mathcal{E}_0(\mathcal{C}_n, \mathcal{D}_n).$$

Hence, using the inequality (34), we have, with probability at least $1 - \delta$,

$$\frac{c}{2}(R_{\text{Bayes}}(\gamma_n^*, Q_n^*) - R_{\text{Bayes}}^*)$$
$$\leq \hat{R}_\phi(\gamma_n^*, Q_n^*) - \hat{R}_\phi(\gamma_n^\dagger, Q_n^\dagger) + 2\mathcal{E}_1(\mathcal{C}_n, \mathcal{D}_n)$$
$$\quad + 2M_n\sqrt{2\ln(2/\delta)/n} + \mathcal{E}_0(\mathcal{C}_n, \mathcal{D}_n)$$
$$\leq 2\mathcal{E}_1(\mathcal{C}_n, \mathcal{D}_n) + \mathcal{E}_0(\mathcal{C}_n, \mathcal{D}_n) + 2M_n\sqrt{2\ln(2/\delta)/n},$$

from which Theorem 2(a) follows.

For part (b), this statement follows by applying (a) with $\delta = 1/n$.

5.3. *Proof of Theorem 3.* One direction of the theorem ("if") is easy. We focus on the other direction. The proof relies on the following technical result.

LEMMA 5. *Given a continuous convex function $f : \mathbb{R}^+ \to \mathbb{R}$, for any $u, v \in \mathbb{R}^+$, define*

$$T_f(u, v) := \left\{ \frac{f^*(\alpha) - f^*(\beta)}{\alpha - \beta} \Big| \alpha \in \partial f(u), \beta \in \partial f(v), \alpha \neq \beta \right\}.$$

*If $f_1 \stackrel{u}{\approx} f_2$, then for any $u, v > 0$, one of the following must be true:*

(1) $T_f(u, v)$ *are nonempty for both $f_1$ and $f_2$, and $T_{f_1}(u, v) = T_{f_2}(u, v)$.*
(2) *Both $f_1$ and $f_2$ are linear in the interval $(u, v)$.*

Now, let us proceed to prove Theorem 3. The convex function $f:[0,\infty) \to \mathbb{R}$ is continuous on $(0,\infty)$ and hence is almost everywhere differentiable on $(0,\infty)$ (see [16]). Note that if function $f$ is differentiable at $u$ and $v$ and $f'(u) \neq f'(v)$, then $T_f(u,v)$ is reduced to a number

$$\frac{uf'(u) - vf'(v) - f(u) + f(v)}{f'(u) - f'(v)} = \frac{f^*(\alpha) - f^*(\beta)}{\alpha - \beta},$$

where $\alpha = f'(u)$, $\beta = f'(v)$, and $f^*$ denotes the conjugate dual of $f$.

Let $v$ be an arbitrary point where both $f_1$ and $f_2$ are differentiable. Let $d_1 = f_1'(v)$, $d_2 = f_2'(v)$. Without loss of generality, we may assume that $f_1(v) = f_2(v) = 0$; if not, we simply consider the functions $f_1(u) - f_1(v)$ and $f_2(u) - f_2(v)$.

Now, for any $u$ where both $f_1$ and $f_2$ are differentiable, applying Lemma 5 for $v$ and $u$, then either $f_1$ and $f_2$ are both linear in $[v,u]$ (or $[u,v]$ if $u<v$), in which case $f_1(u) = cf_2(u)$ for some constant $c$, or the following is true:

$$\frac{uf_1'(u) - f_1(u) - vd_1}{f_1'(u) - d_1} = \frac{uf_2'(u) - f_2(u) - vd_2}{f_2'(u) - d_2}.$$

In either case, we have

$$(uf_1'(u) - f_1(u) - vd_1)(f_2'(u) - d_2) = (uf_2'(u) - f_2(u) - vd_2)(f_1'(u) - d_1).$$

Let $g_1, g_2$ be defined by $f_1(u) = g_1(u) + d_1 u$, $f_2(u) = g_2(u) + d_2 u$. Then, $(ug_1'(u) - g_1(u) - vd_1)g_2'(u) = (ug_2'(u) - g_2(u) - vd_2)g_1'(u)$, implying that $(g_1(u) + vd_1)g_2'(u) = (g_2(u) + vd_2)g_1'(u)$ for any $u$ where $f_1$ and $f_2$ are both differentiable. Since $u$ and $v$ can be chosen almost everywhere, $v$ is chosen so that there does not exist any open interval for $u$ such that $g_2(u) + vd_2 = 0$. It follows that $g_1(u) + vd_1 = c(g_2(u) + vd_2)$ for some constant $c$ and this constant $c$ has to be the same for any $u$ due to the continuity of $f_1$ and $f_2$. Hence, we have $f_1(u) = g_1(u) + d_1 u = cg_2(u) + d_1 u + cvd_2 - vd_1 = cf_2(u) + (d_1 - cd_2)u + cvd_2 - vd_1$. It is now simple to check that $c > 0$ is necessary and sufficient for $I_{f_1}$ and $I_{f_2}$ to have the same monotonicity.

A. *Proof of Lemma 3.* (a) Since $\phi^{-1}(\beta) < +\infty$, we have $\phi(\phi^{-1}(\beta)) = \phi(\inf\{\alpha : \phi(\alpha) \leq \beta\}) \leq \beta$, where the final inequality follows from the lower semi-continuity of $\phi$. If $\phi$ is continuous at $\phi^{-1}(\beta)$, then we have $\phi^{-1}(\beta) = \min\{\alpha : \phi(\alpha) = \beta\}$, in which case we have $\phi(\phi^{-1}(\beta)) = \beta$.

(b) Due to convexity and the inequality $\phi'(0) < 0$, it follows that $\phi$ is a strictly decreasing function in $(-\infty, \alpha^*]$. Furthermore, for all $\beta \in \mathbb{R}$ such that $\phi^{-1}(\beta) < +\infty$, we must have $\phi^{-1}(\beta) \leq \alpha^*$. Therefore, definition 24 and the (decreasing) monotonicity of $\phi$ imply that for any $a, b \in \mathbb{R}$, if $b \geq a \geq \inf \phi$, then $\phi^{-1}(a) \geq \phi^{-1}(b)$, which establishes that $\phi^{-1}$ is a decreasing function. In addition, we have $a \geq \phi^{-1}(b)$ if and only if $\phi(a) \leq b$.

Now, due to the convexity of $\phi$, applying Jensen's inequality for any $0 < \lambda < 1$, we have $\phi(\lambda\phi^{-1}(\beta_1) + (1 - \lambda)\phi^{-1}(\beta_2)) \leq \lambda\phi(\phi^{-1}(\beta_1)) + (1 - \lambda)\phi(\phi^{-1}(\beta_2)) \leq \lambda\beta_1 + (1 - \lambda)\beta_2$. Therefore,

$$\lambda\phi^{-1}(\beta_1) + (1 - \lambda)\phi^{-1}(\beta_2) \geq \phi^{-1}(\lambda\beta_1 + (1 - \lambda)\beta_2),$$

implying the convexity of $\phi^{-1}$.

B. *Proof of Lemma 4.*

PROOF. (a) We first prove the statement for the case of a decreasing function $\phi$. First, if $a \geq b$ and $\phi^{-1}(a) \notin \mathbb{R}$, then $\phi^{-1}(b) \notin \mathbb{R}$; hence, $\Psi(a) = \Psi(b) = +\infty$. If only $\phi^{-1}(b) \notin \mathbb{R}$, then clearly $\Psi(b) \geq \Psi(a)$ [since $\Psi(b) = +\infty$]. If $a \geq b$, and both $\phi^{-1}(\alpha), \phi^{-1}(\beta) \in \mathbb{R}$, then, from the previous lemma, $\phi^{-1}(a) \leq \phi^{-1}(b)$, so that $\phi(-\phi^{-1}(a)) \leq \phi(-\phi^{-1}(b))$, implying that $\Psi$ is a decreasing function.

We next consider the case of a general function $\phi$. For $\beta \in (\beta_1, \beta_2)$, we have $\phi^{-1}(\beta) \in (-\alpha^*, \alpha^*)$, and hence $-\phi^{-1}(\beta) \in (-\alpha^*, \alpha^*)$. Since $\phi$ is strictly decreasing in $(-\infty, \alpha^*]$, then $\phi(-\phi^{-1}(\beta))$ is strictly decreasing in $(\beta_1, \beta_2)$. Finally, when $\beta < \inf \Psi = \phi(\alpha^*)$, $\phi^{-1}(\beta) \notin \mathbb{R}$, so $\Psi(\beta) = +\infty$ by definition.

(b) First of all, assume that $\phi$ is decreasing. By applying Jensen's inequality, for any $0 < \lambda < 1$, we have

$$\begin{aligned}
\lambda\Psi(\gamma_1) &+ (1 - \lambda)\Psi(\gamma_2) \\
&= \lambda\phi(-\phi^{-1}(\gamma_1)) + (1 - \lambda)\phi(-\phi^{-1}(\gamma_2)) \\
&\geq \phi(-\lambda\phi^{-1}(\gamma_1) - (1 - \lambda)\phi^{-1}(\gamma_2)) \quad \text{since } \phi \text{ is convex} \\
&\geq \phi(-\phi^{-1}(\lambda\gamma_1 + (1 - \lambda)\gamma_2)) \\
&= \Psi(\lambda\gamma_1 + (1 - \lambda)\gamma_2),
\end{aligned}$$

where the last inequality is due to the convexity of $\phi^{-1}$ and decreasing $\phi$. Hence, $\Psi$ is a convex function.

In general, the above arguments go through for any $\gamma_1, \gamma_2 \in [\beta_1, \beta_2]$. Since $\Psi(\beta) = +\infty$ for $\beta < \beta_1$, this implies that $\Psi$ is convex in $(-\infty, \beta_2]$.

(c) For any $a \in \mathbb{R}$, from the definition of $\phi^{-1}$ and due to the continuity of $\phi$,

$$\begin{aligned}
\{\beta | \Psi(\beta) = \phi(-\phi^{-1}(\beta)) \leq a\} &= \{\beta | -\phi^{-1}(\beta) \geq \phi^{-1}(a)\} \\
&= \{\beta | \phi^{-1}(\beta) \leq -\phi^{-1}(a)\} \\
&= \{\beta | \beta \geq \phi(-\phi^{-1}(a))\}
\end{aligned}$$

is a closed set. Similarly, $\{\beta \in \mathbb{R} | \Psi(\beta) \geq a\}$ is a closed set. Hence, $\Psi$ is continuous in its domain.

(d) Since $\phi$ is assumed to be classification-calibrated, Lemma 1 implies that $\phi$ is differentiable at 0 and $\phi'(0) < 0$. Since $\phi$ is convex, this implies that $\phi$ is strictly decreasing for $\alpha \leq 0$. As a result, for any $\alpha \geq 0$, let $\beta = \phi(-\alpha)$, then we obtain $\alpha = -\phi^{-1}(\beta)$. Since $\Psi(\beta) = \phi(-\phi^{-1}(\beta))$, we have $\Psi(\beta) = \phi(\alpha)$. Hence, $\Psi(\phi(-\alpha)) = \phi(\alpha)$. Letting $u^* = \phi(0)$, then we have $\Psi(u^*) = u^*$ and $u^* \in (\beta_1, \beta_2)$.

(e) Let $\alpha = \Psi(\beta) = \phi(-\phi^{-1}(\beta)$. Then, from equation (24), $\phi^{-1}(\alpha) \leq -\phi^{-1}(\beta)$. Therefore,

$$\Psi(\Psi(\beta)) = \Psi(\alpha) = \phi(-\phi^{-1}(\alpha)) \leq \phi(\phi^{-1}(\beta)) \leq \beta.$$

We have proved that $\Psi$ is strictly decreasing for $\beta \in (\beta_1, \beta_2)$. As such, $\phi^{-1}(\alpha) = -\phi^{-1}(\beta)$. We also have $\phi(\phi^{-1}(\beta)) = \beta$. It follows that $\Psi(\Psi(\beta)) = \beta$ for all $\beta \in (\beta_1, \beta_2)$.

REMARK. With reference to statement (b), if $\phi$ is not a decreasing function, then the function $\Psi$ need not be convex on the entire real line. For instance, the following loss function generates a function $\Psi$ that is not convex: $\phi(\alpha) = (1-\alpha)^2$ when $\alpha \leq 1$, 0 when $\alpha \in [0,2]$, and $\alpha - 2$ otherwise. Then, we have $\Psi(9) = \phi(2) = 0, \Psi(16) = \phi(3) = 1, \Psi(25/2) = \phi(-1 + 5/\sqrt{2}) = -3 + 5/\sqrt{2} > (\Psi(9) + \Psi(16))/2$. □

C. *Proof of Lemma 5.*

PROOF. Consider a joint distribution $\mathbb{P}(X,Y)$ defined by $\mathbb{P}(Y = -1) = q = 1 - \mathbb{P}(Y = 1)$ and

$$\mathbb{P}(X|Y = -1) \sim \text{Uniform}[0, b] \quad \text{and} \quad \mathbb{P}(X|Y = 1) \sim \text{Uniform}[a, c],$$

where $0 < a < b < c$. Let $\mathcal{Z} = \{1, 2\}$. We assume $Z$ is produced by a deterministic quantizer design $Q$ specified by a threshold $t \in (a, b)$; in particular, we set $Q(z = 1|x) = 1$ when $x \geq t$, and $Q(z = 2|x) = 1$ when $x < t$. Under this quantizer design, we have

$$\mu(1) = (1-q)\frac{t-a}{c-a}; \qquad \mu(2) = (1-q)\frac{c-t}{c-a};$$

$$\pi(1) = q\frac{t}{b}; \qquad \pi(2) = q\frac{b-t}{b}.$$

Therefore, the $f$-divergence between $\mu$ and $\pi$ takes the form

$$I_f(\mu, \pi) = \frac{qt}{b}f\left(\frac{(t-a)b(1-q)}{(c-a)tq}\right) + \frac{q(b-t)}{b}f\left(\frac{(c-t)b(1-q)}{(c-a)(b-t)q}\right).$$

If $f_1 \stackrel{u}{\approx} f_2$, then $I_{f_1}(\mu, \pi)$ and $I_{f_1}(\mu, \pi)$ have the same monotonicity property for any $q \in (0, 1)$, as well, as for any choice of the parameters $q$ and $a < b < c$.

Let $\gamma = \frac{b(1-q)}{(c-a)q}$, which can be chosen arbitrarily positive, and then define the function
$$F(f,t) = tf\left(\frac{(t-a)\gamma}{t}\right) + (b-t)f\left(\frac{(c-t)\gamma}{b-t}\right).$$

Note that the functions $F(f_1, t)$ and $F(f_2, t)$ have the same monotonicity property, for any positive parameters $\gamma$ and $a < b < c$.

We now claim that $F(f, t)$ is a convex function of $t$. Indeed, using convex duality [18], $F(f, t)$ can be expressed as follows:
$$F(f,t) = t \sup_{r \in \mathbb{R}} \left\{ \frac{(t-a)\gamma}{t} r - f^*(r) \right\} + (b-t) \sup_{s \in \mathbb{R}} \left\{ \frac{(c-t)\gamma}{b-t} s - f^*(s) \right\}$$
$$= \sup_{r,s} \{(t-a)r\gamma - tf^*(r) + (c-t)s\gamma - tf^*(s)\},$$

which is a supremum over a linear function of $t$, thereby showing that $F(f, t)$ is convex of $t$.

It follows that both $F(f_1, t)$ and $F(f_2, t)$ are subdifferentiable everywhere in their domains; since they have the same monotonicity property, we must have

(35) $$0 \in \partial F(f_1, t) \Leftrightarrow 0 \in \partial F(f_2, t).$$

It can be verified using subdifferential calculus [8] that
$$\partial F(f,t) = \frac{a\gamma}{t} \partial f\left(\frac{(t-a)\gamma}{t}\right) + f\left(\frac{(t-a)\gamma}{t}\right)$$
$$- f\left(\frac{(c-t)\gamma}{b-t}\right) + \frac{(c-b)\gamma}{b-t} \partial f\left(\frac{(c-t)\gamma}{b-t}\right).$$

Letting $u = \frac{(t-a)\gamma}{t}$, $v = \frac{(c-t)\gamma}{b-t}$, we have

(36a) $\quad 0 \in \partial F(f, t)$

(36b) $\quad \Leftrightarrow \quad 0 \in (\gamma - u)\partial f(u) + f(u) - f(v) + (v - \gamma)\partial f(v)$

(36c) $\quad \Leftrightarrow \quad \exists \alpha \in \partial f(u), \beta \in \partial f(v) \text{ s.t.}$
$$0 = (\gamma - u)\alpha + f(u) - f(v) + (v - \gamma)\beta$$

(36d) $\quad \Leftrightarrow \quad \exists \alpha \in \partial f(u), \beta \in \partial f(v) \text{ s.t.}$
$$\gamma(\alpha - \beta) = u\alpha - f(u) + f(v) - v\beta$$

(36e) $\quad \Leftrightarrow \quad \exists \alpha \in \partial f(u), \beta \in \partial f(v) \text{ s.t. } \gamma(\alpha - \beta) = f^*(\alpha) - f^*(\beta).$

By varying our choice of $q \in (0, 1)$, the number $\gamma$ can take any positive value. Similarly, by choosing different positive values of $a, b, c$ (such that $a < b < c$),

30 X. NGUYEN, M. J. WAINWRIGHT AND M. I. JORDAN

we can ensure that $u$ and $v$ can take on any positive real values such that $u < \gamma < v$. Since equation (35) holds for any $t$, it follows that for any triples $u < \gamma < v$, (36e) holds for $f_1$ if and only if it also holds for $f_2$.

Considering a fixed pair $u < v$, first suppose that the function $f_1$ is linear on the interval $[u, v]$ with a slope $s$. In this case, condition (36e) holds for $f_1$ and any $\gamma$ by choosing $\alpha = \beta = s$, which implies that condition (36e) also holds for $f_2$ for any $\gamma$. Thus, we deduce that $f_2$ is also a linear function on the interval $[u, v]$.

Suppose, on the other hand, that $f_1$ and $f_2$ are both nonlinear in $[u, v]$. Due to the monotonicity of subdifferentials, we have $\partial f_1(u) \cap \partial f_1(v) = \varnothing$ and $\partial f_2(u) \cap \partial f_2(v) = \varnothing$. Consequently, it follows that both $T_{f_1}(u, v)$ and $T_{f_2}(u, v)$ are non-empty. If $\gamma \in T_{f_1}(u, v)$, then condition (36e) holds for $f_1$ for some $\gamma$. Thus, it must also hold for $f_2$ using the same $\gamma$, which implies that $\gamma \in T_{f_2}(u, v)$. The same argument can also be applied with the roles of $f_1$ and $f_2$ reversed, so we conclude that $T_{f_1}(u, v) = T_{f_2}(u, v)$. $\square$

X. Nguyen
Department of Statistical Science
Duke University
Durham, North Carolina 27708
and
Statistical and Applied Mathematical
 Sciences Institute (SAMSI)
Research Triangle Park
Durham, North Carolina 27709
USA
E-mail: xuanlong.nguyen@stat.duke.edu

M. J. Wainwright
M. I. Jordan
Department of Statistics
and
Department of Electrical Engineering
 and Computer Sciences
University of California, Berkeley
Berkeley, California 94720
USA
E-mail: wainwrig@stat.berkeley.edu
 jordan@stat.berkeley.edu